\documentclass[runningheads,envcountsame]{llncs}

\usepackage{graphicx}
\usepackage{xspace}
\usepackage{hyperref}
\usepackage{tikz}
\usepackage{tikzscale}
\usepackage{pgfplots}
\usepackage{boxplot}
\usepackage{amsmath,amssymb}
\usepackage{csquotes}
\usepackage{booktabs}

\pgfplotsset{
  /pgfplots/ybar legend/.style={
    /pgfplots/legend image code/.code={%
      \draw[##1,/tikz/.cd,bar width=3pt,yshift=-0.2em,bar shift=0pt]
    plot coordinates {(2*\pgfplotbarwidth,0.8em)};},
  }
}

\newcommand{\polymake}{\texttt{polymake}\xspace}
\newcommand{\mptopcom}{\texttt{mptopcom}\xspace}
\newcommand{\ppl}{\texttt{PPL}\xspace}

\newcommand{\virosage}{\texttt{Viro.sage}\xspace}

\DeclareMathOperator{\conv}{conv}
\DeclareMathOperator{\pos}{pos}

\newcommand{\RR}{\mathbb R}
\newcommand{\TT}{\mathbb T}
\newcommand{\ZZ}{\mathbb Z}
\newcommand{\CP}{\mathbb{CP}}
\newcommand{\RP}{\mathbb{RP}}
\newcommand{\TP}{\mathbb{TP}}
\newcommand{\1}{{\mathbf 1}}

\newcommand{\cN}{\mathcal N}
\newcommand{\cS}{\mathcal S}
\newcommand{\cP}{\mathcal P}
\newcommand{\cT}{\mathcal T}

\newcommand\SetOf[2]{\left\{#1\,\vphantom{#2}\right|\left.\vphantom{#1}\,#2\right\}}

\newcommand\bigSetOf[2]{\bigl\{#1 \mid #2\bigr\}}

\newcommand\doi[1]{\href{http://dx.doi.org/#1}{\texttt{doi:#1}}}

\usepackage[colorinlistoftodos, bordercolor=orange, backgroundcolor=orange!20, linecolor=orange, textsize=scriptsize]{todonotes}

\title{Real tropical hyperfaces by patchworking \\ in \polymake}
\titlerunning{Patchworking Tropical Hypersurfaces}
\author{Michael Joswig\inst{1,2} \and Paul Vater\inst{2}}
\authorrunning{Joswig and Vater}
\institute{
  TU Berlin, Chair of Discrete Mathematics/Geometry\\
  \and
  MPI MiS Leipzig  \\
  \email{joswig@math.tu-berlin.de},
  \email{vater@mis.mpg.de},\\
}

\begin{document}

\maketitle

% \begin{abstract}
% \end{abstract}


% The main body should describe challenge, achievements and progress in
% mathematical software research, addressing issues such as
% functionality, underlying theories, design, development and applications. 
% at least 4, at most 8 pages

\section{Introduction}

Hilbert's 16th problem asks about topological constraints for real algebraic hypersurfaces in projective space.
In the 1980s Viro developed patchworking as a combinatorial method to construct real algebraic hypersurfaces with unusually large $\ZZ_2$-Betti numbers \cite{Viro:1980,Viro:1983,Viro:1984,Viro:2008}.
A major breakthrough of this idea was Itenberg's refutation of Ragsdale's Conjecture \cite{ItenbergViro:1996}.
Today patchworking is most naturally interpreted within the larger framework of tropical geometry \cite{Tropical+Book}.
In this way patchworking is a combinatorial avenue to real tropical hypersurfaces.

Here we report on a recent implementation of patchworking and real tropical hypersurfaces in \polymake \cite{MPC:polymake}, version 4.1 of April 2020.
The first software for patchworking that we are aware of is the \enquote{Combinatorial Patchworking Tool}~\cite{CombinatorialPatchworkingTool}, which works web-based and is restricted to the planar case.
A second implementation is \virosage \cite{Viro.sage} which is capable of patchworking in arbitrary dimension and degree.
Our implementation has the same scope as \virosage but it is superior in two ways.
First, it naturally ties in with a comprehensive hierarchy of polyhedral objects in \polymake; e.g., this allows for a rich choice of constructions of real tropical hypersurfaces.
Second, our implementation is more efficient.
This is demonstrated by several experiments with curves and surfaces of various degrees.
As a new mathematical contribution we provide a census of Betti numbers of real tropical surfaces.

\paragraph{Acknowledgments.}
We are indebted to Ilia Itenberg, Johannes Rau and Kristin Shaw for valuable comments on a previous version of this text.
Moreover, we are grateful to Lars Kastner and Benjamin Lorenz for helping with the experiments.
M.~Joswig has been supported by DFG (EXC 2046: \enquote{MATH$^+$}, SFB-TRR 195: \enquote{Symbolic Tools in Mathematics and their Application}, and GRK 2434: \enquote{Facets of Complexity}).

\subsection{Tropical hypersurfaces in $\TP^{n-1}$}

Let $f = \bigoplus_{v\in V} c_v \odot x^v \in \TT[x_1,\dots,x_n]$ be a tropical polynomial where $V$ is a finite subset of $\ZZ^n$.
We use the multi-index notation $x^v=x_1^{v_1}\cdots x_n^{v_n}$, and $\TT=\RR\cup\{\infty\}$, $\oplus{=}\min$ and $\odot{=}+$.
The \emph{tropical hypersurface} $\cT(f)$ is the tropical vanishing locus of $f$, i.e., the set of points in $\RR^n$, where the minimum of the evaluation function $x\mapsto f(x)$ is attained at least twice.
Throughout we will assume that $f$ is homogeneous of degree $d$, i.e., for each $v\in V$ we have $v_1+\dots+v_n=d$.
In that case $\cT(f)$ descends to the \emph{tropical projective torus} $\RR^n/\RR\1$, where $\1=(1,\dots,1)$.
The \emph{Newton polytope} of $f$ is $\cN(f)=\conv V$, and the coefficients of $f$ induce a regular subdivision, $\cS(f)$.
The latter is dual to $\cT(f)$.
We refer to \cite{Tropical+Book} and \cite{Triangulations} for further details.

The tropical projective space $\TP^{n-1}=(\TT^n{-}\{\infty\1\})/\RR\1$ compactifies $\RR^n/\RR\1$.
It is naturally stratified into lower dimensional tropical projective tori, marked by those coordinates which are finite.
In this way the pair $(\TP^{n-1},\RR^n/\RR\1)$ is naturally homeomorphic with an $(n{-}1)$-simplex and its interior.
Often we will identify the tropical hypersurface $\cT(f)$ with its compactification in $\TP^{n-1}$.

\subsection{Viro's patchworking}

The following is essentially a condensed version of \cite[\S3.1]{RenaudineauShaw:2018}, with minor variations.
A \emph{sign distribution} $\epsilon \in \ZZ_2^V$ can be \emph{symmetrized} to the function
\[
  s_\epsilon : \ZZ_2^n\to\ZZ_2^V \,,\ s_\epsilon(z)(v) := \epsilon(v) + \langle z,v\rangle \bmod 2 \enspace .
\]
As in \cite{Itenberg:1997} we choose our signs in $\ZZ_2=\{0,1\}$, which corresponds to $\pm1$ via $z\mapsto (-1)^z$.
Further, the elements $z\in\ZZ_2^n$ are in bijection with the $2^n$ orthants of $\RR^n$ via  $z \mapsto \pos\{(-1)^{z_1}e_1,\ldots,(-1)^{z_n}e_n\}$, where $e_1,\dots, e_n$ are the standard basis vectors of $\RR^n$, and $\pos(\cdot)$ denotes the nonnegative hull.
We will use this identification throughout and, consequently, we call $z$ itself an \emph{orthant}.

The tropical hypersurface $\cT(f)$ is a polyhedral complex in $\TP^{n-1}$, and its $k$-dimensional cells are dual to the $(n{-}1{-}k)$-cells of $\cS(f)$.
In particular, each maximal cell $F$ of $\cT(f)$ corresponds to an edge, $V(F)\subset V$, of $\cS(f)$.
We write $\cT_{n-2}$ for the set of maximal cells (which are $(n{-}2)$-dimensional polyhedra) and denote powersets as $\cP(\cdot)$.

\begin{figure}[th]
  \includegraphics[width=.45\textwidth]{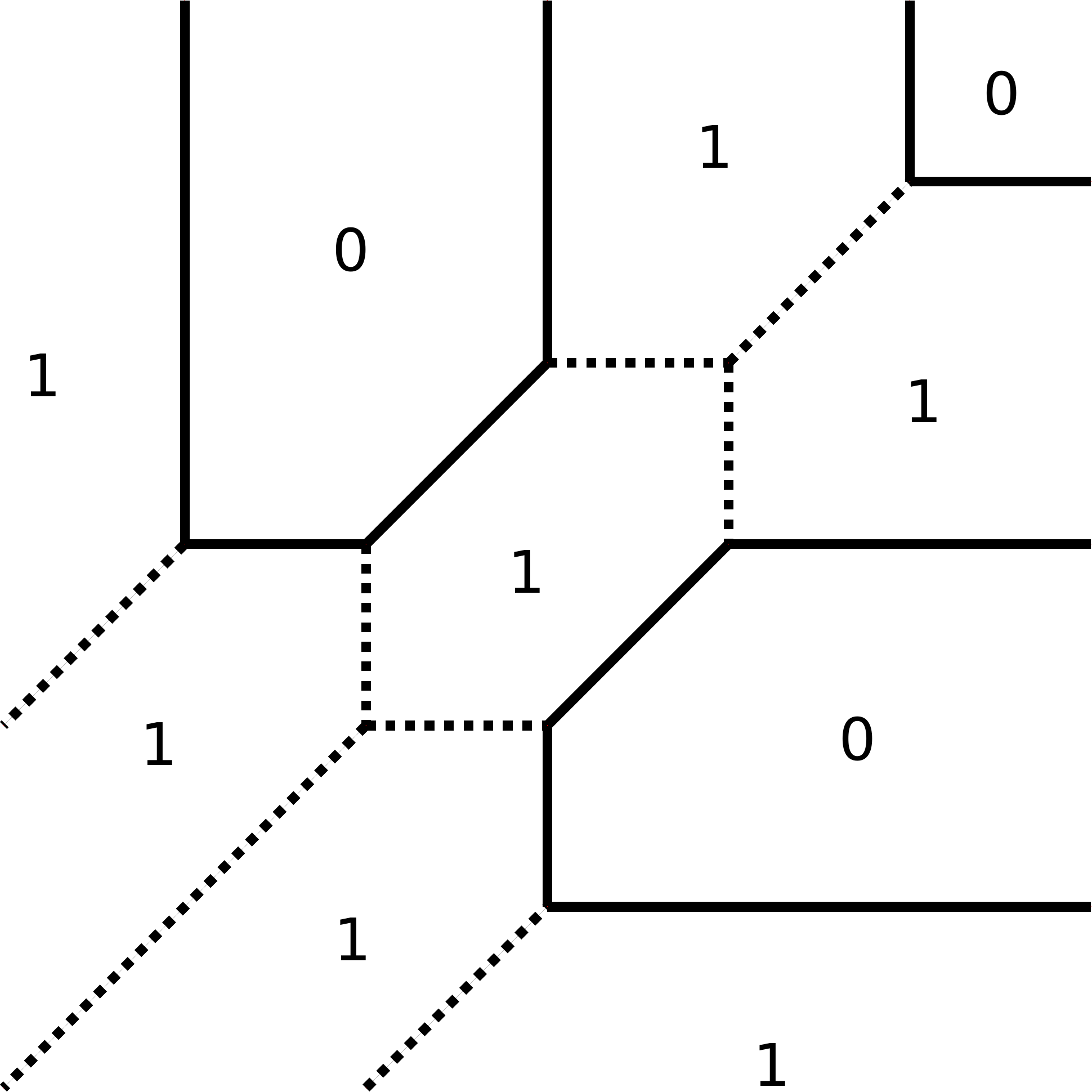} \hfill
  \includegraphics[width=.45\textwidth]{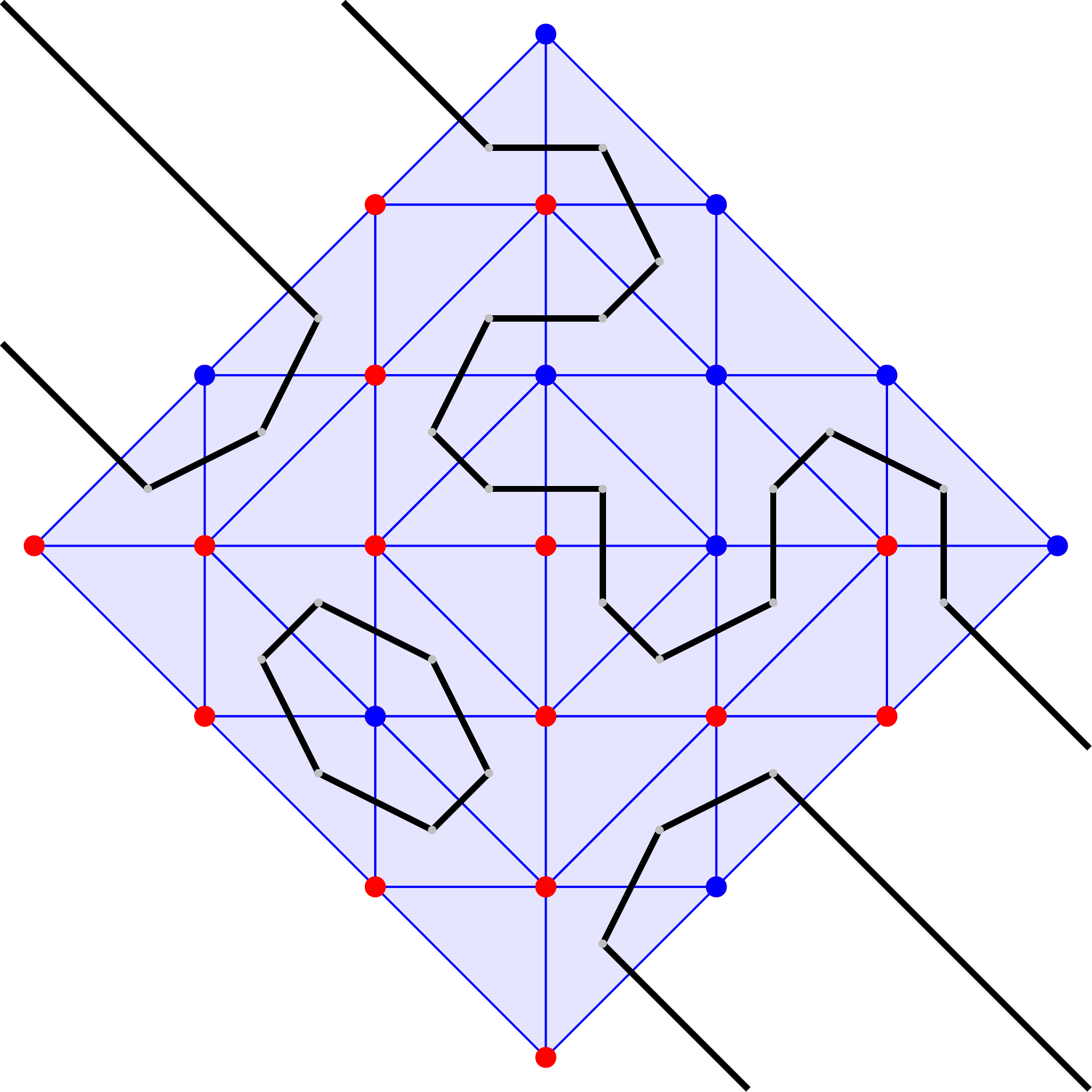}
  \caption{Real tropical elliptic curve (left) and its real part (right)}
  \label{fig:harnack}
\end{figure}

Note that there are no $(n{-}2)$-cells of $\cT(f)$ in the boundary $\TP^{n-1}-\RR^n/\RR\1$.
The \emph{real phase structure} on $\cT(f)$ induced by $\epsilon$ is the map
\[
  \phi_\epsilon : \cT_{n-2}\to\cP(\ZZ_2^n) \,,\, F\mapsto \SetOf{z\in\ZZ_2^n}{s_\epsilon(z)(v)\neq s_\epsilon(z)(w)} \, \text{ for } \{v,w\}=V(F) \enspace .
\]
That is, for each maximal cell $F$ of $\cT(f)$ this describes the set of orthants, in which the symmetrized sign distribution takes distinct values on the two vertices of the dual edge $V(F)$ in $\cS(f)$.
This extends to all cells $G$ of $\cT(f)$ by setting $\phi_\epsilon(G) := \bigcup\phi_\epsilon(F)$, where the union is taken over all maximal cells $F\in\cT_{n-2}$ containing~$G$.
The pair $\cT_\epsilon(f)=(\cT(f),\epsilon)$ is a \emph{real tropical hypersurface}.

Let $\overline{z}$, defined by $\overline{z}_i=1-z_i$, be the \emph{antipode} of $z\in\ZZ_2^n$.
We define an equivalence relation $\sim$ on $\ZZ_2^n\times \TP^{n-1}$, which identifies copies of $\TP^{n-1}$ along common strata, by letting
\[
  (z, x) \sim (z', y) \ :\iff \ x=y \text{ and } \bigl( \overline{z} = z' \text{ or } ( x_i=\infty=y_i \Leftrightarrow z_i = 1 = z_i' ) \bigr) \enspace .
\]
% If we let $\overline{z}$, defined by $\overline{z}_i=1-z_i$, be the \emph{antipode} of $z\in\ZZ_2^n$, we obtain $(z,x)\sim(\overline{z},x)$ whenever $x\not\in\RR^n/\RR\1$ lies in the boundary.
This identifies $\{z\}\times\TP^{n-1}$ and $\{\overline{z}\}\times\TP^{n-1}$ one to one for each $z$.
% and we will consequently use the identification $(\ZZ_2^{n-1}\times \TP^{n-1})/_{\sim} \cong $  $ (\left(\{0\}\times\ZZ_2^{n-1}\right)\times \TP^{n-1})/{\sim} = (\ZZ_2^n\times \TP^{n-1})/{\sim}$.
It follows that the quotient $(\ZZ_2^n\times \TP^{n-1})/{\sim}$ is homeomorphic to the real projective space $\RP^{n-1}$.
Combinatorially that construction can be seen as follows:
the union of the $2^{n}$ simplices $\conv\{(-1)^{z_1}e_1,\ldots,(-1)^{z_{n}}e_{n}\}$, where $z$ ranges over all orthants, gives the boundary of the regular cross polytope $\conv\{\pm e_1,\dots,\pm e_{n}\}$ in~$\RR^{n}$.
Taking the quotient modulo antipodes yields~$\RP^{n-1}$.
% the union of the $2^{n-1}$ simplices $\conv\{0,(-1)^{z_1}e_1,\ldots,(-1)^{z_{n-1}}e_{n-1}\}$, where $z$ ranges over all orthants, gives the regular cross polytope $\conv\{\pm e_1,\dots,\pm e_{n-1}\}$ in~$\RR^{n-1}$.
% Taking the quotient modulo antipodes in the boundary yields~$\RP^{n-1}$.
% Taking first the quotient modulo $\RR\1$ and then modulo antipodes in the boundary yields~$\RP^{n-1}$.

The \emph{real part} of the real tropical hypersurface $\cT_\epsilon(f)=(\cT(f),\epsilon)$, denoted $\RR\cT_\epsilon(f)$, is now defined as the collection of polyhedral complexes in $\ZZ_2^n\times\TP^{n-1}$ consisting of the polyhedra
\[
  \bigSetOf{ \{z\}\times F }{ F\in\cT_{n-2} \text{ and } z \in \phi_\epsilon(F) }
\]
and their faces.
Note that $\{z\}\times F \in \RR\cT_\epsilon(f)$ if and only if $\{\overline{z}\}\times F \in \RR\cT_\epsilon(f)$, and hence we may restrict to the part of $\RR\cT_\epsilon(f)$ in $(\{0\}\times \ZZ_2^{n-1})\times\TP^{n-1}$.

To avoid cumbersome notation and language we call the quotient of $\RR\cT_\epsilon(f)$ by $\sim$ also the \emph{real part} of $\cT_\epsilon(f)$ and use the same symbol, $\RR\cT_\epsilon(f)$.
In this way $\RR\cT_\epsilon(f)$ becomes a piecewise linear hypersurface in $\RP^{n-1}\approx\ZZ_2^n\times\TP^{n-1}/{\sim}$.

% \begin{theorem}[Real structure theorem]
%   [TODO] $\RR\mathcal T_s(f)$ comes as the limit of $\sval$ \ldots
% \end{theorem}

% % We can extend these definitions to projective tropical hypersurfaces:
% % A face $F<\mathcal T(f)|_S$ arises from a unique face $F'<\mathcal T(f)|_\emptyset$, and we simply set $s_\epsilon(F):=s_\epsilon(F')$.

% $\RR\mathcal T_s(f)$ then realizes a topological hypersurface:
% \[
%   \RR\mathcal T_s(f)/_\sim \subset (\ZZ_2^{n+1}\times\TT\RR^n)/_\sim = \RR\PP^n
% \]
% A concrete realization can be constructed in the following way:
% Consider the map $*:\ZZ_2^n\times\RR^n\to\RR^n\,,\, (z,v)\mapsto z*v := ((-1)^{z_1}v_1,\ldots,(-1)^{z_n}v_n)$ [TODO define and use earlier], and the componentwise logarithm $\Log$.

The above construction is relevant for its connection with real algebraic geometry.
To simplify the exposition we now consider a special case:
Setting $\Delta_{n-1}=\conv\{e_1,\dots,e_n\}$, we assume that the set $V=d\cdot\Delta_{n-1}\cap\ZZ^n$ is the set of lattice points in the dilated unit simplex.
This entails that the projective toric variety generated from $V$ is the (complex) projective space $\CP^{n-1}$.
The following result comes in various guises; this version occurs in \cite{Viro:1983} and \cite[Proposition 2.6]{ItenbergShustin:2003}.
\begin{theorem}[Viro's combinatorial patchworking theorem]\label{thm:viro}
  Let $f$ be a homogeneous tropical polynomial of degree $d$ with support $V=d\cdot\Delta_{n-1}\cap\ZZ^n$.
  Then, for each sign distribution $\epsilon\in\ZZ_2^n$, there exists a nonsingular real algebraic hypersurface $X$ in $\CP^{n-1}$, also with Newton polytope $\cN(f)=d\cdot\Delta_{n-1}$, such that
  \[
    (\ZZ_2^n\times\TP^{n-1}/{\sim},\, \RR\cT_\epsilon(f)) \ \text{is $\ZZ_2$-homologous to} \ (\RP^{n-1},\, \RR X) \enspace .
  \]
\end{theorem}

If additionally $\cS(f)$ is \emph{unimodular}, i.e., each simplex has normalized volume one, this is \enquote{primitive patchworking}.
In the primitive case stronger conclusions hold \cite{Viro:1984,RenaudineauShaw:2018}.
The notion \enquote{combinatorial patchworking} refers to the condition $\cN(f)=d\cdot\Delta_{n-1}$.
This is what our implementation supports, for arbitrary degrees and dimensions.
More general results require to carefully take into account the toric geometry of $\cN(f)$.

\begin{example}\label{exmp:harnack}
  With $n=d=3$ we consider the tropical polynomial
  \[
    f \ = \ x^3 \oplus 1x^2y \oplus 1x^2z \oplus 4xy^2 \oplus 3xyz \oplus 4xz^2 \oplus 9y^3 \oplus 7y^2z \oplus 7yz^2 \oplus 9z^3
  \]
  in $\TT[x,y,z]$, where we omit `$\odot$' for improved readability.
  The tropical hypersurface $\cT(f)$ is the tropical elliptic curve in $\RR^3/\RR\1$ in Figure~\ref{fig:harnack} (left).
  The sign distribution $\epsilon=(0, 1, 0, 1, 1, 1, 1, 0, 1, 1)$ yields a real tropical curve with real part in $\ZZ_2^3\times\TP^2/{\sim}$ which has two components; cf.\ Figure\ref{fig:harnack} (right).
  This primitive patchwork corresponds to a classical Harnack curve of degree 3; cf.\ \cite[Sec.~5]{ItenbergViro:1996}.
\end{example}

\section{Betti numbers from combinatorial patchworking}

Our goal is to exhibit a census of Betti numbers of real tropical surfaces in $\ZZ_2^4\times\TP^3/{\sim}$.
Throughout the following let $f$ be a tropical polynomial of degree~$d$ in $n=4$ homogeneous variables; we will assume that $\cS(f)$ is a regular and full triangulation of $V=d\cdot\Delta_3\cap\ZZ^4$.
That is, we focus on combinatorial patchworks.
A triangulation of $V$ is \emph{full} if it uses all points in $V$; a unimodular triangulation is necessarily full.
While the converse holds in the plane, there are many more full triangulations of $d\cdot\Delta_3$ than unimodular ones if $d\geq 3$.
Further, with
\begin{equation}\label{eq:card}
  k \ := \ \frac{1}{6} d^3 + d^2 + \frac{11}{6} d + 1 \enspace,
\end{equation}
which is the cardinality of $V$, we pick a sign vector $\epsilon\in\ZZ_2^k$.
This gives rise to a real algebraic surface $X$ in $\CP^3$ whose real part $\RR X$ is \enquote{near the tropical limit} $\RR\cT_\epsilon(f)$ in the sense of~\cite{RenaudineauShaw:2018}.
Itenberg \cite[Theorems~3.2/3.3]{Itenberg:1997} showed that the Euler characteristic satisfies
\begin{equation}\label{eq:euler}
  \chi(\RR X) \ \geq \ \frac{4d-d^3}{3} \enspace ,
\end{equation}
with equality attained in the primitive/unimodular case.
Moreover, by \cite[Theorem 4.2]{Itenberg:1997},
\begin{equation}\label{eq:betti:1}
  b_1(\RR X) \ \leq \ \frac{2d^3 - 6d^2 + 7d}{3} \enspace,
\end{equation}
where $b_q(\cdot)$ are $\ZZ_2$-Betti numbers; see also \cite{ItenbergShustin:2003:nonfull} for bounds without the fullness assumption.
However, if $\cS(f)$ is even unimodular then, by \cite[Theorem 4.1]{Itenberg:1997},
\begin{equation}\label{eq:betti:0}
  b_0(\RR X) \ \leq \ \binom{d - 1}{3}+1 \enspace.
\end{equation}
See Table~\ref{tab:euler-betti} for explicit numbers in the range which is relevant for our experiments.
The main result of \cite{RenaudineauShaw:2018} furnishes a vast generalization of \eqref{eq:betti:0} to arbitrary dimensions.

\begin{table}[th]\centering
  \caption{Bounds for Euler characteristic and Betti numbers, depending on the degree~$d$.
  The values $k$, $\chi'$, $b_0'$ and $b_1'$ are the right hand sides of \eqref{eq:card}, \eqref{eq:euler}, \eqref{eq:betti:0} and \eqref{eq:betti:1}, respectively.}
  \label{tab:euler-betti}
  \begin{tabular*}{0.75\linewidth}{@{\extracolsep{\fill}}rrrrr@{}}
    \toprule
    $d$ & $k$ & $\chi'$ & $b_0'$ & $b_1'$ \\
    \midrule
    % $2$ & $10$ & $0$     & $1$    & $2$    \\
    $3$ & $20$ & $-5$    & $1$    & $7$    \\
    $4$ & $35$ & $-16$   & $2$    & $20$   \\
    $5$ & $56$ & $-35$   & $5$    & $45$   \\
    $6$ & $84$ & $-64$   & $11$   & $86$   \\
    \bottomrule
  \end{tabular*}
\end{table}

% w z y x
% 0 0 0 3
% 0 0 1 2
% 0 0 2 1
% 0 0 3 0
% 0 1 0 2
% 0 1 1 1
% 0 1 2 0
% 0 2 0 1
% 0 2 1 0
% 0 3 0 0
% 1 0 0 2
% 1 0 1 1
% 1 0 2 0
% 1 1 0 1
% 1 1 1 0
% 1 2 0 0
% 2 0 0 1
% 2 0 1 0
% 2 1 0 0
% 3 0 0 0

% 5 1 1 5 2 0 2 0 0 1 2 0 2 1 1 1 3 3 4 8

% 0 0 0 0 0 0 0 1 1 0 1 1 1 1 1 0 1 1 1 0

\begin{example}\label{exmp:212}
  % We present an example of a patchworked cubic surface with Betti vector $(2,1,2)$.
  The subdivision $\cS(f)$ induced by the tropical polynomial
  \[
    \begin{aligned}
      f \ =& \ 5x^3 \oplus 1x^2y \oplus 1xy^2 \oplus 5y^3 \oplus 2x^2z \oplus 0xyz \oplus 2y^2z \\
      &\oplus 0xz^2 \oplus 0yz^2 \oplus 1z^3 \oplus 2x^2w \oplus 0xyw \oplus 2y^2w \oplus 1xzw \\
      &\oplus 1yzw \oplus 1z^2w \oplus 3xw^2 \oplus 3yw^2 \oplus 4zw^2 \oplus 8w^3
    \end{aligned}
  \]
  is a full triangulation of $3\cdot\Delta_3$ which is not unimodular.
  Its $f$-vector reads $(20,60,64,23)$, and its automorphism group is of order~6.
  The sign distribution
  \[
    \epsilon \ = \ (0, 0, 0, 0, 0, 0, 0, 1, 1, 0, 1, 1, 1, 1, 1, 0, 1, 1, 1, 0)
  \]
  yields a real tropical surface $\RR\cT_\epsilon(f)$ whose real part has Betti vector $(2,1,2)$.
\end{example}

\begin{figure}[ht]
  \centering
  \includegraphics[height=21em]{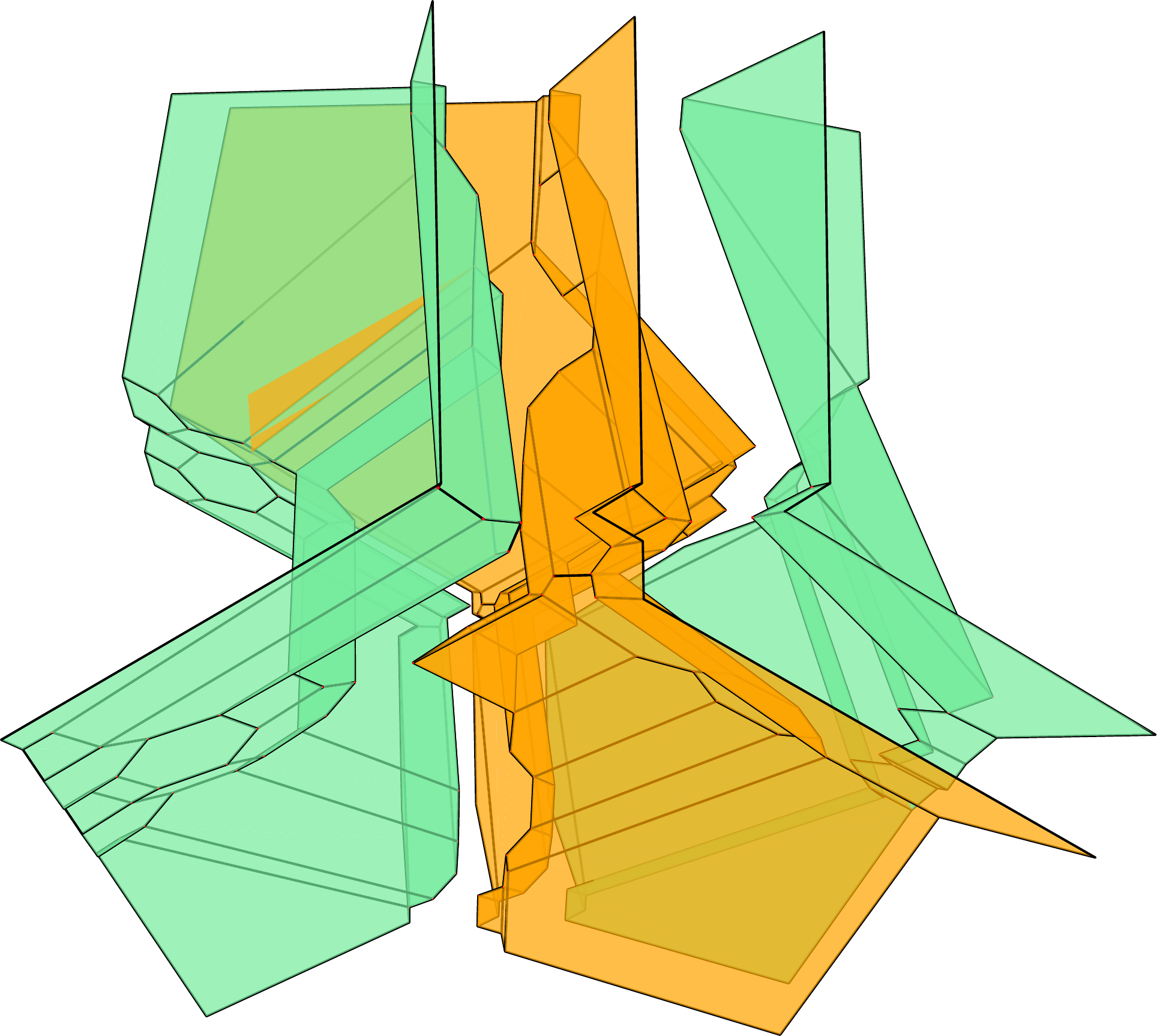}
  \caption{The real part of a cubic surface with Betti vector $(2,1,2)$.
    There are three affine sheets, of which the outer two account for one connected component in $\RP^3$, which is homeomorphic to $\mathbb S^2$; the middle sheet forms a component homeomorphic to $\RP^2$.
  }
  \label{fig:212}
\end{figure}

% $ polymake --script stat.pl 3d3-15001-91301,91303-115001.Betti.xz
% 1 1 1 : 41582 (0.21%)
% 1 3 1 : 745136 (3.73%)
% 1 5 1 : 5708702 (28.54%)
% 1 7 1 : 13503615 (67.52%)
% 2 1 2 : 965 (0.00%)
% #=1000000 total=20000000 avg=20
% ids=2-1000001

% $ polymake --script stat.pl 4d3-1-100000.Betti.xz
% 1 0 1 : 4 (0.00%)
% 1 10 1 : 397129 (19.86%)
% 1 12 1 : 332824 (16.64%)
% 1 14 1 : 162867 (8.14%)
% 1 16 1 : 43575 (2.18%)
% 1 18 1 : 4898 (0.24%)
% 1 2 1 : 134 (0.01%)
% 1 4 1 : 11701 (0.59%)
% 1 6 1 : 92951 (4.65%)
% 1 8 1 : 269536 (13.48%)
% 2 10 2 : 156745 (7.84%)
% 2 12 2 : 204891 (10.24%)
% 2 14 2 : 157081 (7.85%)
% 2 16 2 : 70813 (3.54%)
% 2 18 2 : 17640 (0.88%)
% 2 2 2 : 3 (0.00%)
% 2 20 2 : 1812 (0.09%)
% 2 4 2 : 119 (0.01%)
% 2 6 2 : 10887 (0.54%)
% 2 8 2 : 64390 (3.22%)
% #=100000 total=2000000 avg=20
% ids=1-100000

\begin{figure}[th]
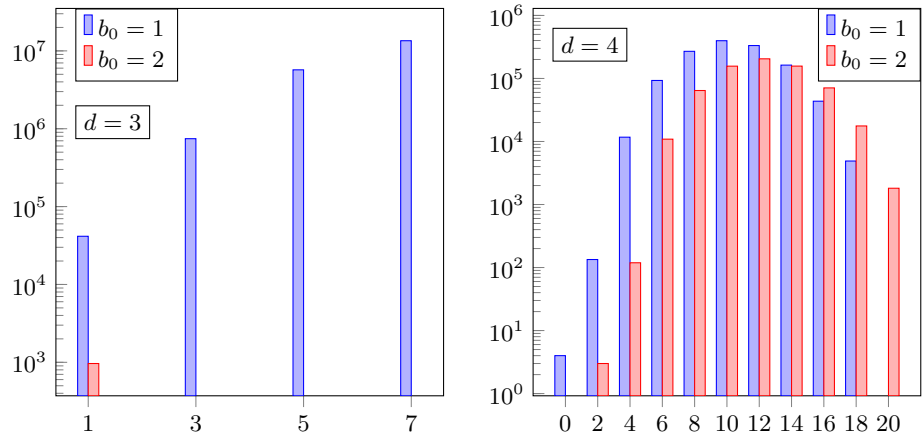

  \centering
  \includegraphics[width=.48\linewidth]{plots/figures/stat.3d3.tikz}\hfill
  \includegraphics[width=.48\linewidth]{plots/figures/stat.4d3.tikz}
  \caption{Distribution of Betti vectors for surfaces of degrees 3 and 4.
    % Values for $b_0=b_2$ are indicated by different colors, values for $b_1$ are marked on the $x$-axis.
    The colors indicate values for $b_0=b_2$, the values on the $x$-axis indicate values for $b_1$.
    For $d=3$ the most frequent vector is $(1,7,1)$ with 67.52\%.
    For $d=4$ it is $(1,10,1)$ with 19.86\%.}
  \label{fig:stat34}
\end{figure}

\subsection{Combinatorial description of the homology}

The polyhedral description of $\RR\cT_\epsilon(f)$ directly gives a combinatorial description of the homology; see also \cite[Proposition 3.17]{RenaudineauShaw:2018}.
The cellular chain modules read
\begin{equation}\label{eq:homology}
  C_q(\RR\cT_\epsilon(f);\ZZ_2) \ = \ \bigoplus_{\sigma \text{ cell of }\cT_\epsilon(f), \dim\sigma = q}\left(\bigoplus_{z\in \phi_\epsilon(\sigma)}\ZZ_2^{\{\sigma\times\{z\}\}}\right)
\end{equation}
and $\partial(\sigma\times\{z\})=\partial(\sigma)\times\{z\}$ defines the boundary maps.
In fact this construction is a special case of a cellular (co-)sheaf \cite{KastnerShawWinz:2017}.
Algorithmically it is beneficial that this does \emph{not} require to geometrically construct $\RR\cT_\epsilon(f)$.

\subsection{A census of Betti numbers of real tropical surfaces}

We used \mptopcom \cite{mptopcom} to compute regular and full triangulations of $d\cdot\Delta_3$ for $3\leq d\leq 6$, which are not necessarily unimodular.
For $d=3$ the total number of such triangulations is known to be $21\,125\,102$ \cite[Table~3]{mptopcom}, up to the natural action of the symmetric group $S_4$.
For higher degrees the corresponding numbers are unknown and probably out of reach for current hard- and software.
Still we can compute some of those triangulations, for each degree.

\begin{figure}[th]
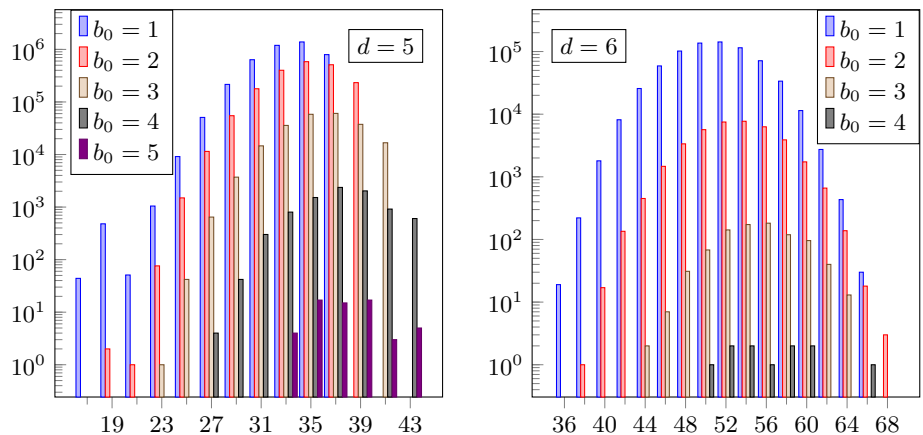

  \centering
  \includegraphics[width=.48\linewidth]{plots/figures/stat.5d3.tikz}\hfill
  \includegraphics[width=.48\linewidth]{plots/figures/stat.6d3.tikz}
  \caption{Distribution of Betti vectors for surfaces of degrees 5 and 6.
    % Values for $b_0=b_2$ are indicated by different colors, values for $b_1$ are marked on the $x$-axis.
    The colors indicate values for $b_0=b_2$, the values on the $x$-axis indicate values for $b_1$.
    For $d=5$ the most frequent vector is $(1,35,1)$ with 21.9\%.
    For $d=6$ it is $(1,52,1)$ with 18.97\%.}
  \label{fig:stat56}
\end{figure}

Our experiments suggest that, in order to see many different Betti vectors $(b_0,b_1,b_0)$, it is preferable to look at many different triangulations.
This is feasible for degrees $3$ and $4$, where we created $1\,000\,000$ and $100\,000$ orbits of triangulations, respectively.
Each of them was equipped with 20 sign distributions which were picked uniformly at random; cf.\ Figure~\ref{fig:stat34}.
For $d=3$ we obtain all values for $b_1$ which are allowed by \eqref{eq:betti:1} if the surface is connected (i.e., $b_0=1$).
Additionally, 965 times we saw the Betti vector $(2,1,2)$; cf.\ Example~\ref{exmp:212}.
In view of \eqref{eq:betti:0} this occurs for non-unimodular triangulations only; all our examples of this kind share the same $f$-vector $(20,60,64,23)$.
For $d=4$ all the possible Betti vectors occur; cf.\ \eqref{eq:euler} and \eqref{eq:betti:1}.

The case of $d=5$ turned out to be surprisingly difficult.
In our standard setup \mptopcom quickly produced about a hundred full and regular triangulations before it stalled.
\mptopcom's algorithm employs a very special search through the flip graph of the point configuration, and it finds all regular triangulations plus some non-regular ones connected by a sequence of flips.
Apparently, most neighbors to our first 100 triangulations of $5\cdot\Delta_3$ are not regular or not full.
As we were interested in exploring many different Betti vectors, we created a second sample of triangulations; to this end we employed a random walk on the flip graph of $5\cdot\Delta_3$.
After eliminating multiples, this gave an additional $13\,000$ regular and full triangulations.
On each of the resulting $13\,100$ triangulations we tried $500$ random sign distributions; cf.\ Figure~\ref{fig:stat56} (left) for the combined statistic.
For $d=6$ we checked $1\,500$ triangulations with $500$ sign distributions each; cf.\ Figure~\ref{fig:stat56} (right).

No matter how hard we try we will only see a tiny fraction of all possible real tropical surfaces of higher degrees.
So the distributions for $d=5$ and $d=6$ may not even be close to the \enquote{truth}.
Yet for $d=5$ we observed $b_1 = 43$, whereas $b_1'=45$; cf.\ Table~\ref{tab:euler-betti}.
We found 61 triangulations of $5\cdot\Delta_3$ with five components, none of which were unimodular.
The maximal number of components in the unimodular case was four.
For $d=6$ our census is way off the theoretical bounds.

\section{Implementation in \polymake}

\polymake is a comprehensive software system for polyhedral geometry and related areas of mathematics \cite{MPC:polymake}.
Mathematical objects like tropical hypersurfaces are determined by their \emph{properties}.
Upon a user query the system directly returns a property (e.g., a tropical polynomial or the dual polyhedral subdivision) if this is known, or it computes it by applying a sequence of \emph{rules}.
Subsequently, the property asked for becomes known, along with any intermediate results.
Throughout the life of such a \emph{big object} the number of properties grows; objects, with their properties, can be saved and loaded again.
The latter is useful, e.g., for processing data on a cluster and examining them on a laptop later.

The computation which is relevant here takes a tropical polynomial $f$ (such that the Newton polytope $\cN(f)$ is a dilated simplex) and a sign distribution $\epsilon$ as input and computes the $\ZZ_2$-Betti numbers of the real part $\RR\cT_\epsilon(f)$ of the real tropical hypersurface $\cT_\epsilon(f)$.
The individual steps are: (1) find the maximal cells of $\cT(f)$ via a dual convex hull computation; (2) compute the Hasse diagram of the entire face lattice of $\cT(f)$; (3) construct the chain complex \eqref{eq:homology} from that Hasse diagram; (4) compute ranks of the boundary matrices mod~2.
Each step is implemented as a separate rule, which makes the code highly modular and reusable.
In particular, the only implementation which is really new is step (3).

We wish to give some details about the first two steps.
Often the dual convex hull computation is the most expensive part.
For this \polymake has interfaces to several algorithms and implementations, the default being \ppl \cite{ppl} which is also used here.
In general, it is difficult to predict which algorithm performs best; see \cite{MPC:polymake} for extensive convex hull experiments.
The computation of the Hasse diagram uses a combinatorial procedure whose complexity is linear in the size of the output, i.e., the total number of cells of the tropical hypersurface; cf.\ \cite{HJS:2019}.

\subsection{Running times}

To compare the running times of \virosage and \polymake for computing the Betti numbers of patchworked hypersurfaces we conducted two experiments, one for Harnack curves and one for surfaces.
All computations were carried out on an AMD Phenom II X6 1090T (3.2GHz, 38528 bmips).

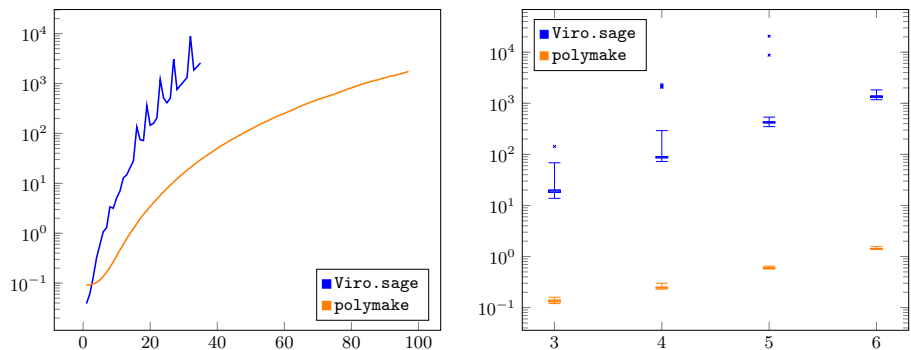
\begin{figure}[th]
  \centering
  \resizebox{0.49\textwidth}{!}{%
    \begin{tikzpicture}
  \begin{axis}[
    xtick pos=left,
    ytick pos=left,
    ymode=log,
    legend pos=south east,
    % thick,
    % legend style={anchor=south west,draw=black,fill=white,legend cell align=left}
    % legend style={legend cell align=left}
    % legend image code/.code={%
    %   \draw[#1, draw=none] (0cm,-0.1cm) -- (0cm,0.1cm);
    % }
    legend cell align={left},
    legend image post style={mark=square*, only marks, mark size=2pt}
    ]
    \addplot[blue,thick] table[mark=none] {plots/figures/time.harnack.avg.virosage.dat};
    \addlegendentry{\virosage}
    \addplot[orange,thick] table[mark=none] {plots/figures/time.harnack.avg.polymake.dat};
    \addlegendentry{\polymake}
  \end{axis}
\end{tikzpicture}
% vim: set ft=tex:
  }%
  \hfill
  \resizebox{0.49\textwidth}{!}{%
    \begin{tikzpicture}
  \begin{axis} [
    box plot width=3pt,
    ymode=log,
    xtick=data,
    % xtick pos=left,
    % ytick pos=left,
    ymode=log,
    % thick,
    legend pos=north west,
    legend cell align={left},
    legend image post style={mark=square*, only marks, mark size=2pt}
    ]

    % \boxplot [forget plot, red] {testdata.dat}
    % \boxplot [blue] {data/testdata2.dat};
    % \addlegendentry{blue}
    \boxplot [blue] {plots/figures/time.xd3.box.virosage.dat};
    \addlegendentry{\virosage}
    \boxplot [orange] {plots/figures/time.xd3.box.polymake.dat};
    \addlegendentry{\polymake}

    \addplot [blue, only marks, mark=x, mark size=1pt] table {plots/figures/time.xd3.box.virosage.outliers.dat};
    % \addplot [black, only marks] table {data/testdata.dat};

    % \boxplot [
    % forget plot,
    % box plot whisker bottom index=1,
    % box plot whisker top index=5,
    % box plot box bottom index=2,
    % box plot box top index=4,
    % box plot median index=3
    % ] {testdata2.dat}

  \end{axis}
\end{tikzpicture}
% vim: set ft=tex:
  }%
  \caption{
    % Time taken (in seconds) to compute Betti numbers ($y$-axis) by degree ($x$-axis).
    Time taken to compute Betti numbers (in seconds).
    % Left: average time for Harnack curves of increasing degree.
    Left: Harnack curves, average time by degree.
    % Right: box plots for various surfaces of degree 3 up to 6.
    Right: various surfaces, boxplots for each degree.
  }
  \label{fig:timesbydeg}
\end{figure}

For the Harnack curves, where we have just one curve per degree (the cubic case is Example~\ref{exmp:harnack}), we repeated the same computation ten times each.
Figure~\ref{fig:timesbydeg} (left) shows the mean running time depending on the degree.
The \virosage code showed a rather wide variety, while the \polymake computations gave almost identical running times for each test.

% left: computation times - harnack curve (i.e., patchworked curve of honeycomb triangulation + harnack signs) of increasing degree, 10 measurements for each degree, on TU server frieda

% right: computation times - deg3 (2000 triangulations), deg4 (1000 t.), deg5 (100 t.), deg6 (75 t.), each for 10 random sign distributions, on TU servers (wandlitz, schwarz, frieda)
% the polymake computations were run on frieda, the viro.sage computations currently come from different servers (of the three mentioned above); at the moment more computations are running on frieda, so that in a few days, when they are finished, we will have all results coming from the same server (frieda)

% user@frieda $ inxi -x -C
% CPU:       Hexa core AMD Phenom II X6 1090T (-MCP-) arch: K10 rev.0 cache: 3072 KB
%            flags: (lm nx sse sse2 sse3 sse4a svm) bmips: 38528
%            clock speeds: max: 3350 MHz 1: 3237 MHz 2: 3230 MHz 3: 3350 MHz 4: 3213 MHz 5: 3218 MHz 6: 3293 MHz

% user@schwarz $ inxi -x -C
% CPU:       Hexa core AMD Phenom II X6 1090T (-MCP-) arch: K10 rev.0 cache: 3072 KB
%            flags: (lm nx sse sse2 sse3 sse4a svm) bmips: 38529
%            clock speeds: max: 3549 MHz 1: 3208 MHz 2: 3549 MHz 3: 3548 MHz 4: 3210 MHz 5: 3545 MHz 6: 3208 MHz

% user@wandlitz $ inxi -x -C
% CPU:       Penta core AMD Phenom II X6 1090T (-MCP-) arch: K10 rev.0 cache: 2560 KB
%            flags: (lm nx sse sse2 sse3 sse4a svm) bmips: 32108
%            clock speeds: max: 3248 MHz 1: 3248 MHz 2: 3213 MHz 3: 3207 MHz 4: 3237 MHz 5: 3241 MHz

The experiment for the surfaces is slightly different in that both the tropical polynomials (and triangulations) and the sign distributions were varied.
% We took the first 2000, 1000, 100, and 75 triangulations for degrees 3, 4, 5, and 6, respectively, and check 10 random sign distributions each.
For degrees 3, 4, 5, and 6 we took the first 2000, 1000, 100, and 75 triangulations (as enumerated by \mptopcom), respectively, and measured the running time for 10 random sign distributions each.
Figure~\ref{fig:timesbydeg} (right) shows a box plot for each degree.
The boxes indicate the 2nd and 3rd quartiles, the whiskers mark the minimum and maximum time measurements, excluding outliers (i.e., measurements whose ratio to the median is either bigger than $4$, or smaller than $0.25$), which are marked separately.
Again \virosage exhibits a much greater variety of running times than \polymake.

\section{Conclusion}

We have shown that our new implementation is capable of determining the $\ZZ_2$-Betti numbers of a patchworked surface of moderate degree within a few seconds.
This allows for providing a rich census.

One major reason for \polymake being faster than \virosage \cite{Viro.sage} is that we avoid the explicit construction of a simplicial complex model of $\RR\cT_\epsilon(f)$.
Moreover, \polymake computes $\ZZ_2$ Betti numbers directly, while \virosage goes through a standard homology computation with integer coefficients.
\polymake provides geometric realizations (and integral homology), too, but this is unnecessary here.

\bibliographystyle{splncs04}
\bibliography{patchworking}

\end{document}